# On Maxwell fluid with relaxation time and viscosity depending on the pressure

Satish Karra

*Texas A&M University*
*Department of Mechanical Engineering*
*3123 TAMU, College Station TX 77843-3123*
*United States of America*

Vít Průša[1,*]

*Faculty of Mathematics and Physics*
*Charles University in Prague*
*Sokolovská 83*
*Praha 8 – Karlín, CZ 186 75*
*Czech Republic*

K. R. Rajagopal[**]

*Texas A&M University*
*Department of Mechanical Engineering*
*3123 TAMU, College Station TX 77843-3123*
*United States of America*

**Abstract**

We study a variant of the well known Maxwell model for viscoelastic fluids, namely we consider the Maxwell fluid with viscosity and relaxation time depending on the pressure. Such a model is relevant for example in modelling behaviour of some polymers and geomaterials. Although it is experimentally known that the material moduli of some viscoelastic fluids can depend on the pressure, most of the studies concerning the motion of viscoelastic fluids do not take such effects into account despite their possible practical significance in technological applications. Using a generalized Maxwell model with pressure dependent material moduli we solve a simple boundary value problem and we demonstrate interesting non-classical features exhibited by the model.

## 1. Introduction

The Maxwell fluid model was originally developed by Maxwell [1] to describe the elastic and viscous response of air. Nowadays, it is however, frequently used to model the response of various viscoelastic fluids ranging from polymers—see for example Ferry [2]—to the Earth's mantle—see for example Cathles [3]. In the present paper, we study an important generalization of the original model due to Maxwell,

---

[*]Current address: Texas A&M University, Department of Mechanical Engineering, 3123 TAMU, College Station, TX 77843-3123, United States of America.

[**]Corresponding author.

*Email addresses:* satkarra@tamu.edu (Satish Karra), prusv@karlin.mff.cuni.cz (Vít Průša), krajagopal@tamu.edu (K. R. Rajagopal)

[1]Vít Průša thanks the Nečas Center for Mathematical Modeling (project LC06052 finaced by the MŠMT of the Czech republic) for its support.

namely we consider a model with pressure dependent material moduli. While the model with a pressure dependent material moduli has important technological ramifications, little is known about the qualitative or quantitative features related to the model. In fact, there is no careful analytical study with regard to mathematical questions concerning existence and uniqueness of solutions for such fluids. Even within the context of solutions to initial-boundary value problems, there is no systematic study when all the material moduli are pressure dependent. Given the possible usefulness of this model, it is surprising that there are no such studies and the analysis carried out here addresses this lacuna.

Bridgman [4] in his pioneering experiments in high pressure physics reported (for many organic fluids) a significant dependence of the viscosity on the pressure and articulated the need for the dependence of the material moduli on the pressure. Models of fluids with pressure dependent viscosity are nowadays frequently used to describe the behaviour of fluids in many applications, for example in lubrication theory, see Neale [5] and Gwynllyw et al. [6]. Extending to viscoelastic fluids, we can ask whether the material moduli of viscoelastic fluid models also exhibit dependence on pressure, and what effects can be described if we use pressure dependent material moduli. Viscoelastic material moduli depending on the pressure has been reported for polymers, see for example Singh and Nolle [7], McKinney and Belcher [8] and the literature stemming from these papers, as well as for geomaterials, see for example Weertman et al. [9], Ivins et al. [10] and Sahaphol and Miura [11].

The question on pressure dependent viscosity and/or relaxation time is especially interesting, for example, with respect to applications in geophysics, since the material of the Earth's mantle is subject to a wide range of pressures. If we consider a material stratified due to the influence of the gravitational force, then because of pressure dependent material moduli, we would be dealing with a body with material moduli depending on the vertical coordinate. A similar situation occurs in geophysical applications and usually it is assumed that the body (in this case the Earth's mantle) is composed of a number of layers of materials with constant material moduli—a paradigm introduced in papers by McConnel [12, 13]. The approach based on material moduli dependence on the pressure and consequently (in the case of a stratified material) on the vertical coordinate provides an alternative approach to the problem. In such an approach, the material moduli continuously vary with the depth, in contrast to McConnel [13]. Although geophysicists and polymer engineers are aware of the possibility of pressure dependent material moduli, this dependence is invariably ignored in studies concerning the dynamics of these materials. If we take into account that geologists often try to refine their models for the viscosity by considering advanced non-Newtonian models, for example a power-law type viscosity—see for example Weertman et al. [9], Wu and Wang [14]—one wonders why important physical phenomenon such as material moduli depending on pressure is not considered.

In the present paper we would like to explore potential benefits of models based on the assumption that the material moduli are pressure dependent. We consider a simple boundary value problem for an incompressible Maxwell fluid with the relaxation time and/or viscosity depending on the pressure to illustrate the consequences of the material moduli depending upon the pressure. The constitutive model therefore has the following structure:

$$\mathbb{T} = -\pi\mathbb{I} + \mathbb{S}, \tag{1.1a}$$

$$\mathbb{S} + \lambda(\pi)\overset{\triangledown}{\mathbb{S}} = 2\mu(\pi)\mathbb{D}, \tag{1.1b}$$

where $\mathbb{T}$ denotes the Cauchy stress tensor, $\pi$ is the pressure[2].

$$\mu(\pi) = \mu_0 \left(1 + \beta\pi\right), \tag{1.2a}$$
$$\mu(\pi) = \mu_0 e^{\beta\pi}, \tag{1.2b}$$
$$\lambda(\pi) = \lambda_0 \left(1 + \gamma\pi\right), \tag{1.2c}$$
$$\lambda(\pi) = \lambda_0 e^{\gamma\pi}, \tag{1.2d}$$

where $\lambda_0$, $\mu_0$, $\beta$ and $\gamma$ are given constants. Viscosity in the form (1.2a) and (1.2b) is used in studies dealing with viscous fluids with pressure dependent viscosity, see for example Hron et al. [17] and Srinivasan and Rajagopal [18]; models (1.2c) and (1.2d) are counterparts of (1.2a) and (1.2b) with respect to modulus $\lambda$. From the theoretical point of view, all the models that are considered fit into the thermodynamical framework developed by Rajagopal and Srinivasa [19] (see also Rajagopal and Srinivasa [15]).

The boundary value problem we are going to solve using the model (1.1) is a classical one; it is a variant of the second Stokes' problem introduced by Stokes [20]. The fluid is confined between two parallel plates, the bottom (top) plate is at rest and the top (bottom) plate is moving with time-periodic velocity $\mathbf{V} = V \cos(\omega t)\mathbf{e}_{\hat{x}}$. The pressure on the bottom plate is fixed to $\pi_0$ and a specific body force—the gravitational force—is acting on the fluid, the gravitational acceleration is denoted as $g$. (See Figure 1.) We call the problem a variant of the second Stokes' problem, since Stokes [20] considered the problem in a half-space above an oscillating plate and without the presence of the gravitational force. A setting identical to ours was previously studied by Srinivasan and Rajagopal [18] in their paper concerning viscous fluids with pressure dependent viscosities.

There is a further need for the careful assessment of the model by considering initial-boundary value problems that have technological relevance and also those that can be used to correlate the model with experiments. From the point of geophysical applications another simple problem that is worth studying using the model (1.1) is a variant of problem introduced by Haskell [21] (response of a fluid on the removal of load), and isostatic adjustment in incompressible, nonrotating and self-gravitating spherical planet, a problem introduced by Love [22]—see also discussion in Wolf [23].

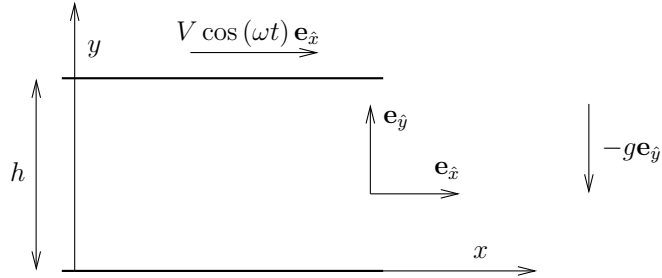

Figure 1: Problem geometry.

---

[2]In this paper, we shall refer to the Lagrange multiplier $\pi$ as pressure. For the model (1.1), the Lagrange multiplier is not the same as the mean normal stress, and if one defines the mechanical pressure as being the mean normal stress, then $\pi$ is not the mechanical pressure. For further discussion see Rajagopal and Srinivasa [15], Rajagopal [16].

## 2. Governing equations

We now turn our attention to solving the boundary value problem. Equations governing the motion of the material are

$$\rho \frac{d\mathbf{v}}{dt} = \operatorname{div} \mathbb{T} + \rho \mathbf{b}, \tag{2.1a}$$

$$\operatorname{div} \mathbf{v} = 0, \tag{2.1b}$$

where $\mathbf{v}$ denotes the velocity, $\rho$ is the density, $\mathbf{b} = -g\mathbf{e}_{\hat{y}}$ is the specific body force and $\mathbb{T}$ is given by constitutive model (1.1). No-slip boundary conditions on the top and bottom plate read—for the oscillating top plate—$\mathbf{v}|_{y=0} = 0$ and $\mathbf{v}|_{y=h} = V\cos(\omega t)\mathbf{e}_{\hat{x}}$, and vice versa for the oscillating bottom plate[3]. Since we are dealing with a simple geometry, the system (2.1) can be greatly simplified under the assumption of parallel flow.

### 2.1. Governing equations for parallel flow

Let us assume that the velocity field and the pressure field have the form $\mathbf{v} = v^{\hat{x}}(y,t)\mathbf{e}_{\hat{x}}$, $\pi = \pi(y)$ respectively, and let us suppose that $\mathbb{S} = \mathbb{S}(y,t)$. It follows that $[\nabla \mathbf{v}]\mathbf{v} = 0$, $[\nabla \mathbb{S}]\mathbf{v} = 0$ and

$$\mathbb{L} = \begin{bmatrix} 0 & \frac{\partial v^{\hat{x}}}{\partial y} \\ 0 & 0 \end{bmatrix}, \quad \mathbb{D} = \frac{1}{2}\begin{bmatrix} 0 & \frac{\partial v^{\hat{x}}}{\partial y} \\ \frac{\partial v^{\hat{x}}}{\partial y} & 0 \end{bmatrix}, \quad \overset{\triangledown}{\mathbb{S}} = \begin{bmatrix} \frac{\partial S_{\hat{x}\hat{x}}}{\partial t} - 2\frac{\partial v^{\hat{x}}}{\partial y}S_{\hat{x}\hat{y}} & \frac{\partial S_{\hat{x}\hat{y}}}{\partial t} - \frac{\partial v^{\hat{x}}}{\partial y}S_{\hat{y}\hat{y}} \\ \frac{\partial S_{\hat{x}\hat{y}}}{\partial t} - \frac{\partial v^{\hat{x}}}{\partial y}S_{\hat{y}\hat{y}} & \frac{\partial S_{\hat{y}\hat{y}}}{\partial t} \end{bmatrix}.$$

The balance of linear momentum (2.1a) reduces to

$$\rho \frac{\partial v^{\hat{x}}}{\partial t} = \frac{\partial S_{\hat{x}\hat{y}}}{\partial y}, \tag{2.2a}$$

$$0 = -\frac{\partial \pi}{\partial y} + \frac{\partial S_{\hat{y}\hat{y}}}{\partial y} - \rho g, \tag{2.2b}$$

and (1.1b) reduces to

$$S_{\hat{x}\hat{x}} + \lambda(\pi)\frac{\partial S_{\hat{x}\hat{x}}}{\partial t} - 2\lambda(\pi)\frac{\partial v^{\hat{x}}}{\partial y}S_{\hat{x}\hat{y}} = 0, \tag{2.2c}$$

$$S_{\hat{x}\hat{y}} + \lambda(\pi)\frac{\partial S_{\hat{x}\hat{y}}}{\partial t} - \lambda(\pi)\frac{\partial v^{\hat{x}}}{\partial y}S_{\hat{y}\hat{y}} = \mu(\pi)\frac{\partial v^{\hat{x}}}{\partial y}, \tag{2.2d}$$

$$S_{\hat{y}\hat{y}} + \lambda(\pi)\frac{\partial S_{\hat{y}\hat{y}}}{\partial t} = 0. \tag{2.2e}$$

Obviously, if we are interested in a time-periodic solution, we need to fix $S_{\hat{y}\hat{y}} = 0$. If $S_{\hat{y}\hat{y}} = 0$, we can then solve (2.2b) to get a formula for the pressure

$$\pi = -\rho g y + \pi_0. \tag{2.3}$$

---

[3]For the sake of brevity we will, however, describe a solution procedure only for the top plate oscillating, and we will only report the results for the case when the bottom plate is oscillating.

Finally, the system reduces to

$$\rho \frac{\partial v^{\hat{x}}}{\partial t} = \frac{\partial S_{\hat{x}\hat{y}}}{\partial y}, \qquad (2.4a)$$

$$S_{\hat{x}\hat{x}} + \lambda(\pi)\frac{\partial S_{\hat{x}\hat{x}}}{\partial t} - 2\lambda(\pi)\frac{\partial v^{\hat{x}}}{\partial y}S_{\hat{x}\hat{y}} = 0, \qquad (2.4b)$$

$$S_{\hat{x}\hat{y}} + \lambda(\pi)\frac{\partial S_{\hat{x}\hat{y}}}{\partial t} = \mu(\pi)\frac{\partial v^{\hat{x}}}{\partial y}, \qquad (2.4c)$$

where $\lambda(\pi)$ and $\mu(\pi)$ are, by virtue of (2.3), known functions of $y$. Equations (2.4a) and (2.4c) are linear equations for $v^{\hat{x}}$ and $S_{\hat{x}\hat{y}}$ and these equations can be solved independently on (2.4b). The system has to be completed by supplying boundary conditions, in our case (oscillating top plate) the boundary conditions for the velocity read

$$\left. v^{\hat{x}} \right|_{y=h} = V\cos(\omega t), \qquad \left. v^{\hat{x}} \right|_{y=0} = 0. \qquad (2.5)$$

### 2.2. Dimensionless governing equations for parallel flow

Let us develop a dimensionless version of the governing equations (2.4) before we proceed to solving the problem. Let $h$ be the characteristic length and $V$ the characteristic velocity, then system (2.4) reads

$$\frac{\partial v^{\hat{x}\star}}{\partial t^\star} = \frac{1}{\mathrm{Re}}\frac{\partial S^\star_{\hat{x}\hat{y}}}{\partial y^\star}, \qquad (2.6a)$$

$$S^\star_{\hat{x}\hat{x}} + \mathrm{We}\,\lambda^\star(\pi^\star)\left(\frac{\partial S^\star_{\hat{x}\hat{x}}}{\partial t^\star} - 2\frac{\partial v^{\hat{x}\star}}{\partial y^\star}S^\star_{\hat{x}\hat{y}}\right) = 0, \qquad (2.6b)$$

$$S^\star_{\hat{x}\hat{y}} + \mathrm{We}\,\lambda^\star(\pi^\star)\frac{\partial S^\star_{\hat{x}\hat{y}}}{\partial t^\star} = \mu^\star(\pi^\star)\frac{\partial v^{\hat{x}\star}}{\partial y^\star}, \qquad (2.6c)$$

where $\mathrm{Re} = \frac{\rho V h}{\mu_0}$ is the Reynolds number (the ratio of inertial forces to viscous forces), $\mathrm{We} = \frac{\lambda_0 V}{h}$ is the Weissenberg number (the ratio of the relaxation time of the material and a characteristic time scale), and the star denotes dimensionless variables $v^{\hat{x}\star} = \frac{v^{\hat{x}}}{V}$, $S^\star_{\hat{x}\hat{x}} = \frac{hS_{\hat{x}\hat{x}}}{\mu_0 V}$, $S^\star_{\hat{x}\hat{y}} = \frac{hS_{\hat{x}\hat{y}}}{\mu_0 V}$, $y^\star = \frac{y}{h}$, $t^\star = \frac{Vt}{h}$. Functions $\lambda^\star(\pi^\star)$ and $\mu^\star(\pi^\star)$ are given by formulae

$$\mu^\star(\pi^\star) = (1 + \beta^\star \pi^\star), \qquad (2.7a)$$

$$\mu^\star(\pi^\star) = \mathrm{e}^{\beta^\star \pi^\star}, \qquad (2.7b)$$

$$\lambda^\star(\pi^\star) = (1 + \gamma^\star \pi^\star), \qquad (2.7c)$$

$$\lambda^\star(\pi^\star) = \mathrm{e}^{\gamma^\star \pi^\star}. \qquad (2.7d)$$

where $\beta^\star = \beta\pi_0$, $\gamma^\star = \gamma\pi_0$ and

$$\pi^\star = -\Pi y^\star + 1, \qquad (2.8)$$

where $\Pi = \frac{h\rho g}{\pi_0}$ is a dimensionless version of (2.3). Clearly, the dimensionless version of boundary conditions (2.5) reads

$$\left. v^{\hat{x}\star} \right|_{y^\star=0} = 0, \qquad \left. v^{\hat{x}\star} \right|_{y^\star=1} = \cos(\omega^\star t^\star). \qquad (2.9)$$

Hereafter, we will use only dimensionless variables and we will therefore omit the star denoting the dimensionless variables. Furthermore, in what follows we will freely vary all the dimensionless parameters that are present in the equations. The parameter values and their mutual relationship do not necessarily correspond to a realistic setting; the main aim of the present study is to illustrate the trends that correspond to particular variations of the parameters.

# 3. Solution to the governing equations

Due to the linearity of the governing equations and the nature of the boundary conditions (2.5), we seek a solution of the form

$$v^{\hat{x}}(y, t) = \tilde{v}^{\hat{x}}(y)e^{i\omega t}, \quad S_{\hat{x}\hat{y}}(y, t) = \tilde{S}_{\hat{x}\hat{y}}(y)e^{i\omega t}, \tag{3.1}$$

where $\tilde{v}^{\hat{x}}(y)$ is a complex function, solution to (2.6) is then the real part of $\tilde{v}^{\hat{x}}(y)e^{i\omega t}$ and $\tilde{S}_{\hat{x}\hat{y}}(y)e^{i\omega t}$ respectively. Substituting (3.1) to (2.6) leads to system

$$i\omega \tilde{v}^{\hat{x}} = \frac{1}{\text{Re}} \frac{d\tilde{S}_{\hat{x}\hat{y}}}{dy}, \tag{3.2a}$$

$$(1 + i\omega \text{We}\lambda(\pi))\tilde{S}_{\hat{x}\hat{y}} = \mu(\pi)\frac{d\tilde{v}^{\hat{x}}}{dy}, \tag{3.2b}$$

that can be reduced to a single ordinary differential equation for the velocity component $\tilde{v}^{\hat{x}}$,

$$i\omega \tilde{v}^{\hat{x}} = \frac{1}{\text{Re}} \frac{d}{dy}\left(\frac{\mu(\pi)}{1 + i\omega \text{We}\lambda(\pi)} \frac{d\tilde{v}^{\hat{x}}}{dy}\right), \tag{3.3}$$

that has to be solved subject to boundary conditions

$$\Re\left(\tilde{v}^{\hat{x}}e^{i\omega t}\right)\big|_{y=0} = 0, \quad \Re\left(\tilde{v}^{\hat{x}}e^{i\omega t}\right)\big|_{y=1} = \cos(\omega t), \tag{3.4}$$

where $\Re$ denotes the real part of the corresponding expression. The boundary conditions must be satisfied for all $t$. Equation (3.3) with boundary conditions (3.4) can be rewritten as a system of ordinary differential equations for the real and imaginary part of $\tilde{v}^{\hat{x}}$. Let $\tilde{v}^{\hat{x}}_{\text{Re}}$ and $\tilde{v}^{\hat{x}}_{\text{Im}}$ denote the real and imaginary part of $\tilde{v}^{\hat{x}}$, then (3.3) reads as follows

$$-\omega \text{Re}\tilde{v}^{\hat{x}}_{\text{Im}} = \frac{d}{dy}\left[\frac{\mu(\pi)}{1 + (\omega \text{We}\lambda(\pi))^2} \frac{d\tilde{v}^{\hat{x}}_{\text{Re}}}{dy} + \frac{\mu(\pi)\omega \text{We}\lambda(\pi)}{1 + (\omega \text{We}\lambda(\pi))^2} \frac{d\tilde{v}^{\hat{x}}_{\text{Im}}}{dy}\right], \tag{3.5a}$$

$$\omega \text{Re}\tilde{v}^{\hat{x}}_{\text{Re}} = \frac{d}{dy}\left[-\frac{\mu(\pi)\omega \text{We}\lambda(\pi)}{1 + (\omega \text{We}\lambda(\pi))^2} \frac{d\tilde{v}^{\hat{x}}_{\text{Re}}}{dy} + \frac{\mu(\pi)}{1 + (\omega \text{We}\lambda(\pi))^2} \frac{d\tilde{v}^{\hat{x}}_{\text{Im}}}{dy}\right], \tag{3.5b}$$

Using (3.1), the required velocity is given by

$$v^{\hat{x}}_{\text{Re}} = \tilde{v}^{\hat{x}}_{\text{Re}} \cos(\omega t) - \tilde{v}^{\hat{x}}_{\text{Im}} \sin(\omega t), \tag{3.6}$$

and the boundary conditions (3.4) reduce to

$$\tilde{v}^{\hat{x}}_{\text{Re}}\big|_{y=0} = 0, \quad \tilde{v}^{\hat{x}}_{\text{Re}}\big|_{y=1} = 1, \quad \tilde{v}^{\hat{x}}_{\text{Im}}\big|_{y=0} = 0, \quad \tilde{v}^{\hat{x}}_{\text{Im}}\big|_{y=1} = 0. \tag{3.7}$$

In addition, the vorticity only has a component is the z-direction whose value is given by

$$(\text{rot }\mathbf{v})^{\hat{z}} = -\frac{dv^{\hat{x}}_{\text{Re}}}{dy} = -\frac{d\tilde{v}^{\hat{x}}_{\text{Re}}}{dy}\cos(\omega t) + \frac{d\tilde{v}^{\hat{x}}_{\text{Im}}}{dy}\sin(\omega t). \tag{3.8}$$

## 3.1. Analytical solution to the governing equations

Let us now solve (3.3) subject to boundary conditions (3.7). Although we will not be able to get a closed analytical formulae for the solution, we can, in most cases, give analytical formulae wherein certain integration constants have to be found numerically.

### 3.1.1. Case $\lambda(\pi) = \lambda_0$, $\mu(\pi) = \mu_0$

In the classical case (3.3) reduces to

$$i\omega \tilde{v}^{\hat{x}} = \frac{1}{\mathrm{Re}} \frac{1}{1 + i\omega \mathrm{We}} \frac{\mathrm{d}^2 \tilde{v}^{\hat{x}}}{\mathrm{d}y^2}, \tag{3.9}$$

and the solution reads

$$\tilde{v}^{\hat{x}}(y) = C_1 \mathrm{e}^{\mathrm{i}^{\frac{1}{2}} By} + C_2 \mathrm{e}^{-\mathrm{i}^{\frac{1}{2}} By}, \tag{3.10}$$

where $B = \sqrt{\omega \mathrm{Re}\,(1 + \mathrm{i}\omega \mathrm{We})}$, and $C_1$ and $C_2$ are arbitrary complex constants. In this case it is in principle possible to apply the boundary conditions (3.7) and obtain values of $C_1$ and $C_2$ as combinations of elementary functions of parameters Re, We and $\omega$. (A solution to the classical second Stokes' problem in the infinite domain was given by Fetecau and Fetecau [24].) Also, note that in this case the velocity field is not influenced by the gravitational force; however, the gravitational force affects the pressure field. In the cases that follow, the gravitational force affects both velocity and pressure fields.

### 3.1.2. Case $\lambda(\pi) = \lambda_0 (1 + \gamma\pi)$, $\mu(\pi) = \mu_0$

If $\lambda(\pi)$ is given by (2.7c), then using (2.8) we can substitute for the pressure and get the relaxation time $\lambda(\pi)$ as a function of $y$, $\lambda(y) = 1 + \gamma - \gamma \Pi y$. Denoting $\hat{y} = 1 + \mathrm{i}\omega \mathrm{We}\lambda(\pi)$, and using $\hat{y}$ as a new independent variable, $\tilde{v}^{\hat{x}}(y) = \hat{v}^{\hat{x}}(\hat{y})$, we see that (3.3) can be rewritten as

$$\hat{y}\frac{\mathrm{d}^2 \hat{v}^{\hat{x}}}{\mathrm{d}\hat{y}^2} - \frac{\mathrm{d}\hat{v}^{\hat{x}}}{\mathrm{d}\hat{y}} + \mathrm{i}\frac{\mathrm{Re}}{\omega\,(\mathrm{We}\gamma\Pi)^2}\hat{y}^2 \hat{v}^{\hat{x}} = 0. \tag{3.11}$$

This equation is of form $x^2 y'' + axy' + (bx^n + c) y = 0$, $b \neq 0$, $n \neq 0$, for which the solution (see Polyanin and Zaitsev [25]) has the form $y = x^{\frac{1-a}{2}}\left(C_1 \mathrm{J}_\nu\left(\frac{2}{n}\sqrt{b}x^{\frac{n}{2}}\right) + C_2 \mathrm{Y}_\nu\left(\frac{2}{n}\sqrt{b}x^{\frac{n}{2}}\right)\right)$, where $\nu = \frac{1}{n}\sqrt{(1-a)^2 - 4c}$, $C_1$ and $C_2$ are arbitrary constants, and $\mathrm{J}_\nu(z)$ and $\mathrm{Y}_\nu(z)$ denote Bessel functions of the first and second kind. Solution to (3.11) therefore reads

$$\tilde{v}^{\hat{x}}(y) = \hat{y}\left(C_1 \mathrm{J}_{\frac{2}{3}}\left(\frac{2}{3}\mathrm{i}^{\frac{1}{2}} B\hat{y}^{\frac{3}{2}}\right) + C_2 \mathrm{Y}_{\frac{2}{3}}\left(\frac{2}{3}\mathrm{i}^{\frac{1}{2}} B\hat{y}^{\frac{3}{2}}\right)\right), \tag{3.12}$$

where $B = \frac{1}{\mathrm{We}\gamma\Pi}\sqrt{\frac{\mathrm{Re}}{\omega}}$, $\hat{y} = 1 + \mathrm{i}\omega \mathrm{We}\,(1 + \gamma - \gamma\Pi y)$, and $C_1$ and $C_2$ are arbitrary complex constants that can be fixed by applying the boundary conditions (3.4). Since there is no general formula for splitting a Bessel function of a general complex argument into its imaginary and real parts (see for example Watson [26]), the system of linear algebraic equations arising from the application of the boundary conditions contains factors—the Bessel functions evaluated at $y = 0$ and $y = 1$—that must be found numerically. An algorithm for numerical computation of Bessel functions of the complex argument is discussed for example by Amos [27]. Formula (3.12) is worth of considering in the case when it is necessary to have a kind of analytical expression for the solution. If one however needs only a numerical solution, it is simpler to directly solve the system (3.5), see Section 3.2. We will face the same situation in the remaining cases.

### 3.1.3. Case $\lambda(\pi) = \lambda_0$, $\mu(\pi) = \mu_0 (1 + \beta\pi)$

If $\mu(\pi)$ is given by (2.7a), then using (2.8) we can substitute for the pressure and get the viscosity $\mu(\pi)$ as a function of $y$, $\mu(y) = 1 + \beta - \beta\Pi y$. Denoting $\hat{y} = (1 + \beta - \beta\Pi y)^{\frac{1}{2}}$, and using $\hat{y}$ as a new independent variable, $\tilde{v}^{\hat{x}}(y) = \hat{v}^{\hat{x}}(\hat{y})$, we see that (3.3) can be rewritten as

$$\hat{y}^2 \frac{\mathrm{d}^2 \hat{v}^{\hat{x}}}{\mathrm{d}\hat{y}^2} + \hat{y}\frac{\mathrm{d}\hat{v}^{\hat{x}}}{\mathrm{d}\hat{y}} - \frac{4\mathrm{i}\omega\mathrm{Re}\,(1 + \mathrm{i}\omega\mathrm{We})}{(\beta\Pi)^2}\hat{y}^2 \hat{v}^{\hat{x}} = 0. \tag{3.13}$$

This equation is again of form $x^2 y'' + axy' + (bx^n + c) y = 0$, $b \neq 0$, $n \neq 0$. Solution to (3.13) therefore reads

$$\tilde{v}^{\hat{x}}(y) = \hat{y} \left( C_1 \mathrm{J}_0 \left( \mathrm{i}^{\frac{3}{2}} B \hat{y} \right) + C_2 \mathrm{Y}_0 \left( \mathrm{i}^{\frac{3}{2}} B \hat{y} \right) \right), \tag{3.14}$$

where $B = \frac{2}{\beta \Pi} \sqrt{\omega \mathrm{Re}(1 + \mathrm{i}\omega \mathrm{We})}$, $\hat{y} = (1 + \beta - \beta \Pi y)^{\frac{1}{2}}$, and $C_1$ and $C_2$ are arbitrary complex constants that can be fixed by applying the boundary conditions (3.4).

### 3.1.4. Case $\lambda(\pi) = \lambda_0 \mathrm{e}^{\gamma \pi}$, $\mu(\pi) = \mu_0$

If $\lambda(\pi)$ is given by (2.7d), then using (2.8) we can substitute for the pressure and get the relaxation time $\lambda(\pi)$ as a function of $y$, $\lambda(y) = \mathrm{e}^{\gamma(1 - \Pi y)}$. Denoting $\hat{y} = \mathrm{i}\omega \mathrm{We}\, \mathrm{e}^{\gamma(1 - \Pi y)}$, and using $\hat{y}$ as a new independent variable, $\tilde{v}^{\hat{x}}(y) = \hat{v}^{\hat{x}}(\hat{y})$, we see that (3.3) can be rewritten as

$$\hat{y}^2 (1 + \hat{y}) \frac{\mathrm{d}^2 \hat{v}^{\hat{x}}}{\mathrm{d}\hat{y}^2} - \hat{y}^2 \frac{\mathrm{d}\hat{v}^{\hat{x}}}{\mathrm{d}\hat{y}} - \mathrm{i} B (1 + \hat{y}) \hat{v}^{\hat{x}} = 0, \tag{3.15}$$

where $B = \frac{\mathrm{Re}}{\omega (\mathrm{We} \gamma \Pi)}$. Equation (3.15) to the best of our knowledge, does not have an analytical solution. It is however possible to find a solution in terms of power series, but it would bring an additional complication to those discussed in Section 3.1.2, and we will therefore not give a solution in terms of power series.

### 3.1.5. Case $\lambda(\pi) = \lambda_0$, $\mu(\pi) = \mu_0 \mathrm{e}^{\beta \pi}$

If $\mu(\pi)$ is given by (2.7b), then using (2.8) we can substitute for the pressure and get the viscosity $\mu(\pi)$ as a function of $y$, $\mu(\pi) = \mathrm{e}^{\beta - \beta \Pi y}$. Denoting $\hat{y} = (\mu(\pi))^{-\frac{1}{2}} = \mathrm{e}^{-\frac{1}{2}(\beta - \beta \Pi y)}$, and using $\hat{y}$ as a new independent variable, $\tilde{v}^{\hat{x}}(y) = \hat{v}^{\hat{x}}(\hat{y})$, we see that (3.3) can be rewritten as

$$\hat{y}^2 \frac{\mathrm{d}^2 \hat{v}^{\hat{x}}}{\mathrm{d}\hat{y}^2} - \hat{y} \frac{\mathrm{d}\hat{v}^{\hat{x}}}{\mathrm{d}\hat{y}} - \frac{4\mathrm{i}\omega \mathrm{Re}(1 + \mathrm{i}\omega \mathrm{We})}{(\beta \Pi)^2} \hat{y}^2 \hat{v}^{\hat{x}} = 0. \tag{3.16}$$

This equation is again of form $x^2 y'' + axy' + (bx^n + c) y = 0$, $b \neq 0$, $n \neq 0$. Solution to (3.16) therefore reads

$$\tilde{v}^{\hat{x}}(y) = \hat{y} \left( C_1 \mathrm{J}_1 \left( \mathrm{i}^{\frac{3}{2}} B \hat{y} \right) + C_2 \mathrm{Y}_1 \left( \mathrm{i}^{\frac{3}{2}} B \hat{y} \right) \right), \tag{3.17}$$

where $B = \frac{2}{\beta \Pi} \sqrt{\omega \mathrm{Re}(1 + \mathrm{i}\omega \mathrm{We})}$, $\hat{y} = \mathrm{e}^{-\frac{1}{2}(\beta - \beta \Pi y)}$, and $C_1$ and $C_2$ are arbitrary complex constants that can be fixed by applying the boundary conditions (3.4).

## 3.2. Numerical solution

For the sake of convenience, we shall denote $\frac{\mu(\pi)}{1 + (\omega \mathrm{We} \lambda(\pi))^2}$ by $\mathcal{A}$ and $\frac{\mu(\pi) \omega \mathrm{We} \lambda(\pi)}{1 + (\omega \mathrm{We} \lambda(\pi))^2}$ by $\mathcal{B}$, and the terms in the square brackets of the right hand side of (3.5a) and (3.5b) by $\mathcal{F}$ and $\mathcal{G}$ respectively. Then, (3.5) can be re-written as a system of four first-order ordinary differential equations as follows

$$\frac{\mathrm{d}\mathcal{F}}{\mathrm{d}y} = -\omega \mathrm{Re} \tilde{v}^{\hat{x}}_{\mathrm{Im}}, \tag{3.18a}$$

$$\frac{\mathrm{d}\mathcal{G}}{\mathrm{d}y} = \omega \mathrm{Re} \tilde{v}^{\hat{x}}_{\mathrm{Re}}, \tag{3.18b}$$

$$\frac{\mathrm{d}\tilde{v}^{\hat{x}}_{\mathrm{Re}}}{\mathrm{d}y} = \frac{\mathcal{AF} - \mathcal{BG}}{\mathcal{A}^2 + \mathcal{B}^2}, \tag{3.18c}$$

$$\frac{\mathrm{d}\tilde{v}^{\hat{x}}_{\mathrm{Im}}}{\mathrm{d}y} = \frac{\mathcal{AG} + \mathcal{BF}}{\mathcal{A}^2 + \mathcal{B}^2}. \tag{3.18d}$$

The system of ordinary differential equations (3.18) along with the boundary conditions (3.7) was solved using the solver `bvp4c` in MATLAB, for various cases of $\mu(\pi)$, $\lambda(\pi)$.

Plots of the solutions for various parameter values are given in Figures 2–14. Parameter values were chosen to demonstrate trends induced by their change and do not necessarily correspond to parameter values encountered in real situations. Figures 2–11 display solution for oscillating top plate, while Figures 12–14 display the solution when the bottom plate is oscillating.

Let us now consider the case when the top plate is oscillating. Snapshots of velocity profiles for various models are shown in Figure 2. Models (1.2a), (1.2c) and (1.2d), for given parameter values, lead to very similar velocity profiles, see Figures 2a, 2c–2e. Model (1.2b), however, for the given parameter values, shows a substantial departure from the classical model (constant viscosity and constant relaxation time), compare Figure 2a and Figure 2b. The difference between the velocity profiles is clearly visible in Figure 3a where we compare, at a particular time instant, the classical model and model (1.2c) for various values of $\gamma$. It is obvious that model (1.2c) allows the oscillations induced by the top plate to propagate up to the bottom plate. This holds for the other non-classical models as well, see Figure 3b and Figure 4, but model (1.2b) leads to strongest oscillations near the bottom plate. Moreover, model (1.2b) leads to the creation of layers where the velocity rapidly changes its direction with respect to the vertical coordinate.

The fact that (1.2b) leads to the propagation of oscillations up to the bottom plate is also well documented by plots of vorticity given in Figures 7, 8, 9 and Figure 11. (The vorticity has, in the case of parallel flow, only one nonzero component, rot $\mathbf{v} = -\frac{\partial v^{\hat{x}}}{\partial y}\mathbf{e}_{\hat{z}}$.) Heavily oscillating layers are generated in the fluid (1.2b) especially for high Weissenberg number (We), Reynolds number (Re) and, of course, $\gamma$; see Figure 5 and Figure 6 for the comparison of velocity profiles, and Figure 10 and Figure 11 for the comparison of vorticities. Obviously, with sufficiently large We, Re and $\gamma$, the model leads to a motion where the vorticity is concentrated near the bottom plate, although it is the top plate that forces the material to move, see Figure 11c.

Similar effects can be observed in the case where the bottom plate is oscillating and the top plate is at rest, see Figures 12–14. Model (1.2b) again leads to the creation of highly oscillating layers near the bottom plate. On the other hand, model (1.2d) with increasing $\beta$, leads to nearly vanishing velocity gradient near the bottom plate and steep velocity gradient near the top plate, see Figure 14a; one can say that the material develops a boundary layer close to the boundary opposite to the boundary inducing the motion (compare cases $\beta = 0$ and $\beta = 8$)!

## 4. Conclusion

We have studied the time dependent flow of a generalized Maxwell model that can be used to describe viscoelastic materials in which the material moduli depend on the pressure. Such models are relevant in modelling behaviour of polymers and geomaterials. Using the generalized model we have solved a simple boundary value problem—a variant of the well known Stokes' problem for the flow induced by an oscillating plate. The solutions were found numerically, and for some particular cases we have also derived analytical expressions for the solution. We have shown that the pressure dependent material moduli reduce to material moduli that are continuously dependent on the vertical coordinate, and that the non-constant material moduli, compared to the classical Maxwell model, have substantial impact on the dynamical behaviour of the fluid. The departures from the behaviour predicted by the classical model are the most distinctive for the model with constant viscosity and relaxation time depending exponentially on the pressure. If the latter model is used instead of the classical one for the case when the top plate is oscillating, then the oscillations induced by the moving top plate can propagate up to the bottom plate, and the vorticity can reach its maximum near the bottom plate; these effects however cannot be captured by the classical Maxwell model.

Although the fact that material moduli for viscoelastic fluids can depend on the pressure is known, it is rarely used in dynamical considerations, even if one expects that the arising pressure to be significantly different at different parts of the body. The present paper provides a simple illustration of the effects that

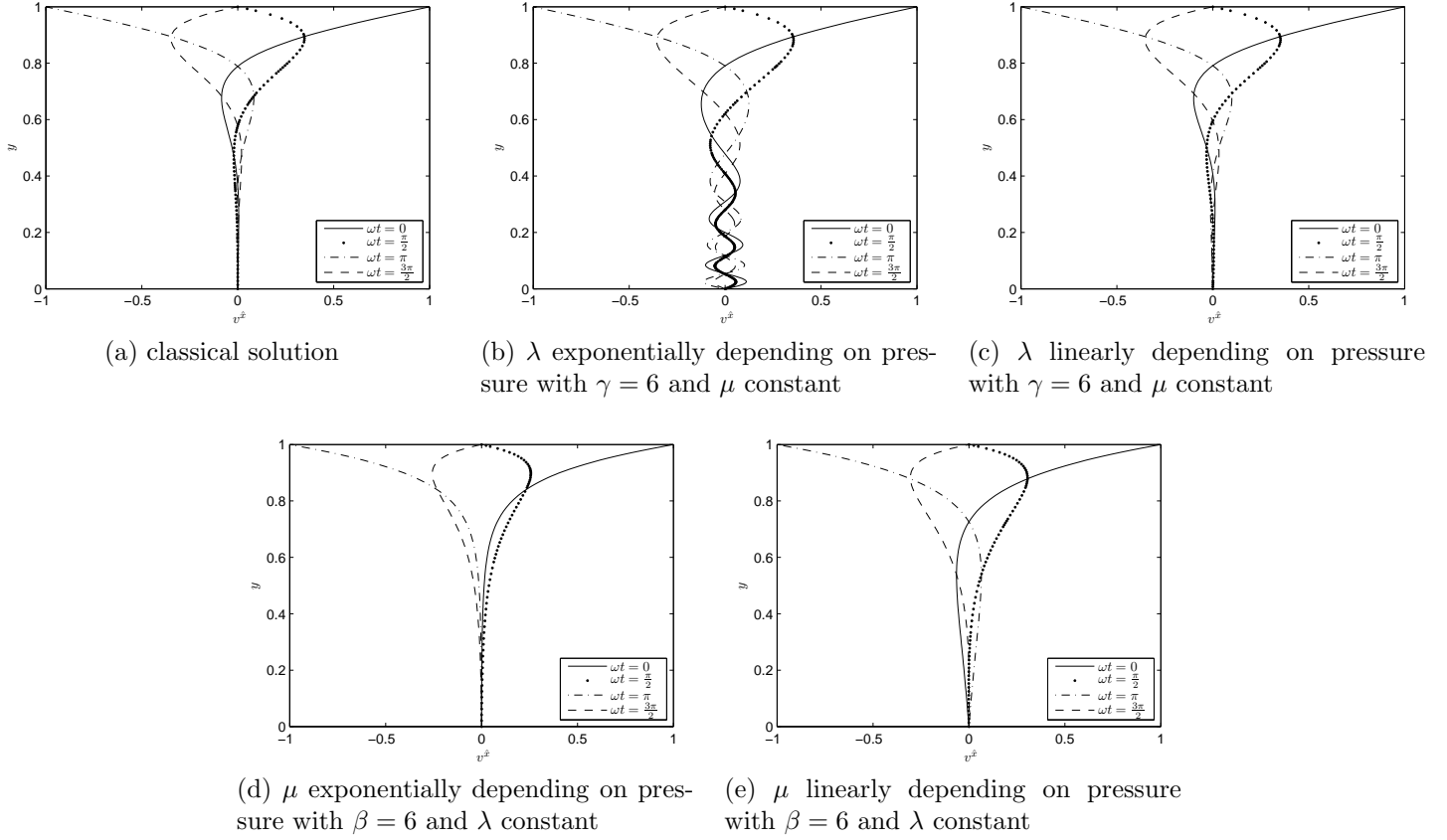

Figure 2: Velocity profiles for various cases of relaxation time ($\lambda$) and viscosity ($\mu$) depending on pressure at different times, Re = 100, We = 0.1, $\omega = 1$, $\Pi = 1$.

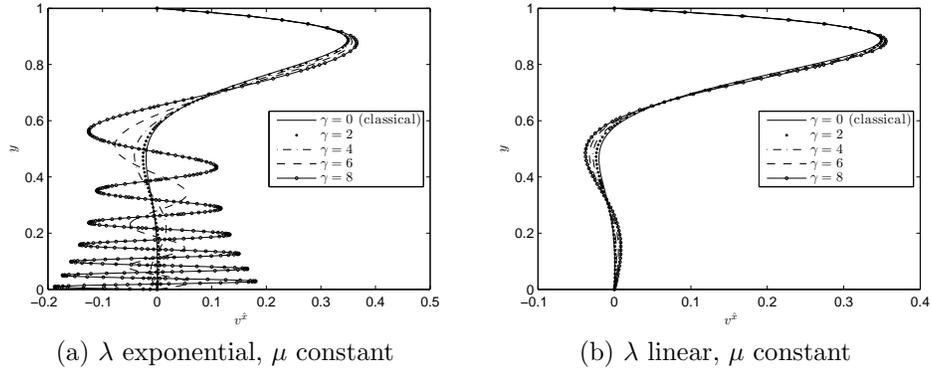

Figure 3: Velocity profiles at $\omega t = \frac{\pi}{2}$ for various $\gamma$ values when the relaxation time ($\lambda$) depends on pressure exponentially and linearly. Viscosity ($\mu$) is kept constant in both cases with Re = 100, We = 0.1, $\omega = 1$, $\Pi = 1$.

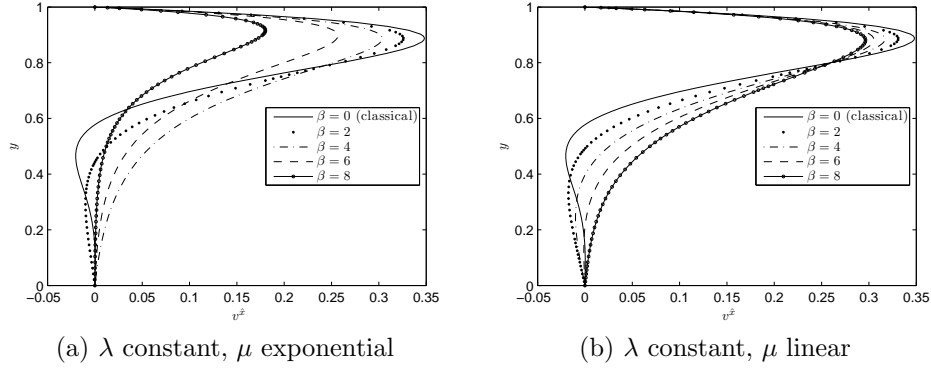

(a) $\lambda$ constant, $\mu$ exponential

(b) $\lambda$ constant, $\mu$ linear

Figure 4: Velocity profiles at $\omega t = \frac{\pi}{2}$ for various $\beta$ values when the viscosity ($\mu$) depends on pressure exponentially and linearly. Relaxation time ($\lambda$) is kept constant in both cases with Re = 100, We = 0.1, $\omega = 1$, $\Pi = 1$.

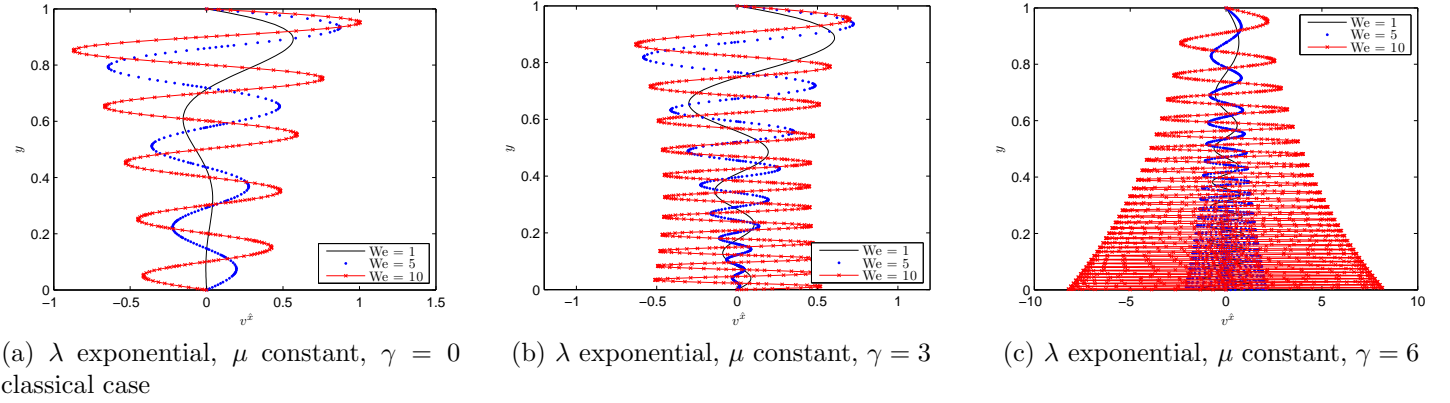

(a) $\lambda$ exponential, $\mu$ constant, $\gamma = 0$ classical case

(b) $\lambda$ exponential, $\mu$ constant, $\gamma = 3$

(c) $\lambda$ exponential, $\mu$ constant, $\gamma = 6$

Figure 5: Velocity profiles at $\omega t = \frac{\pi}{2}$ for various We values for the classical case as well as when the relaxation time ($\lambda$) depends on pressure exponentially. Viscosity ($\mu$) is kept constant in all cases with Re = 100, $\omega = 1$, $\Pi = 1$.

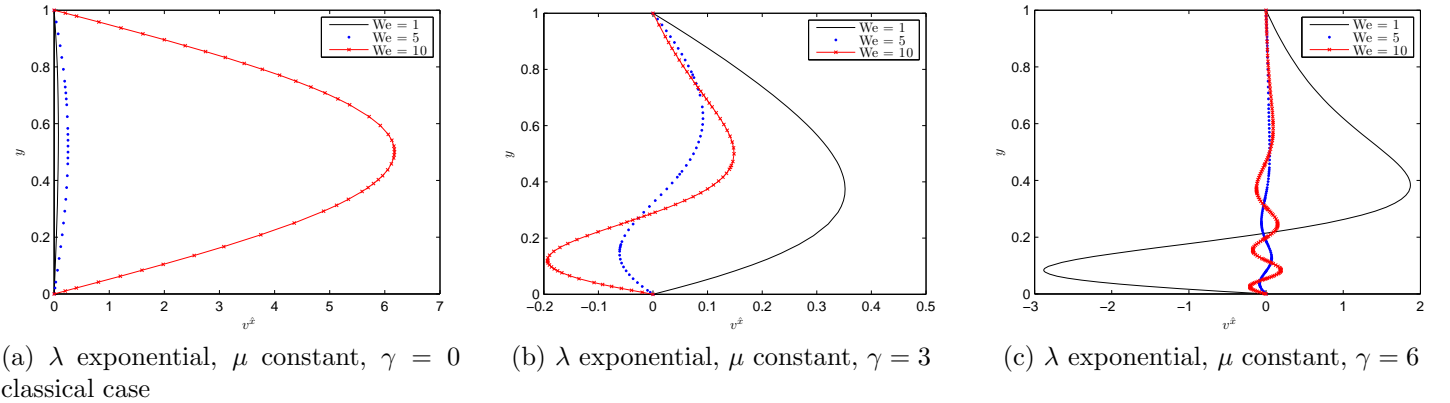

(a) $\lambda$ exponential, $\mu$ constant, $\gamma = 0$ classical case

(b) $\lambda$ exponential, $\mu$ constant, $\gamma = 3$

(c) $\lambda$ exponential, $\mu$ constant, $\gamma = 6$

Figure 6: Velocity profiles at $\omega t = \frac{\pi}{2}$ for various We values for the classical case as well as when the relaxation time ($\lambda$) depends on pressure exponentially. Viscosity ($\mu$) is kept constant in all cases with Re = 1, $\omega = 1$, $\Pi = 1$.

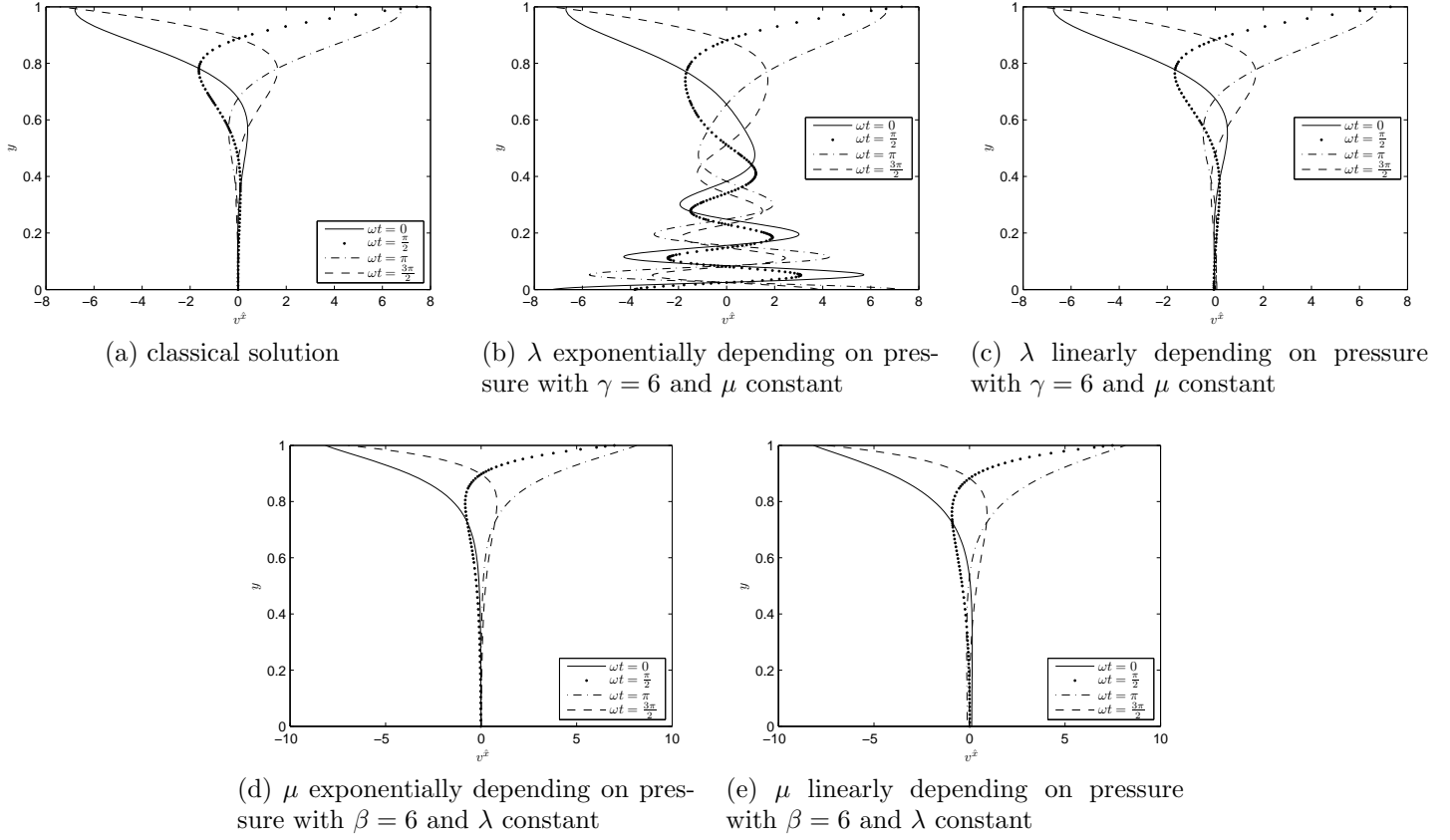

Figure 7: Snapshots of vorticity profiles for various cases of relaxation time ($\lambda$) and viscosity ($\mu$) depending on pressure, Re = 100, We = 0.1, $\omega = 1$, $\Pi = 1$.

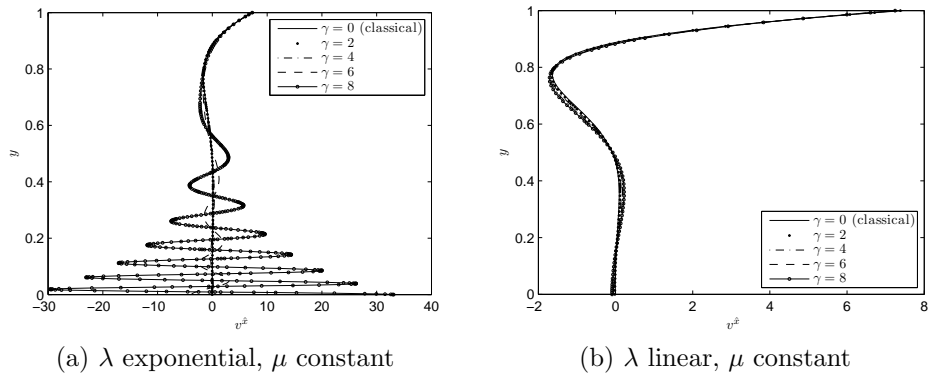

Figure 8: Vorticity profiles at $\omega t = \frac{\pi}{2}$ for various $\gamma$ values when the relaxation time ($\lambda$) depends on pressure exponentially and linearly. Viscosity ($\mu$) is kept constant in both cases with Re = 100, We = 0.1, $\omega = 1$, $\Pi = 1$.

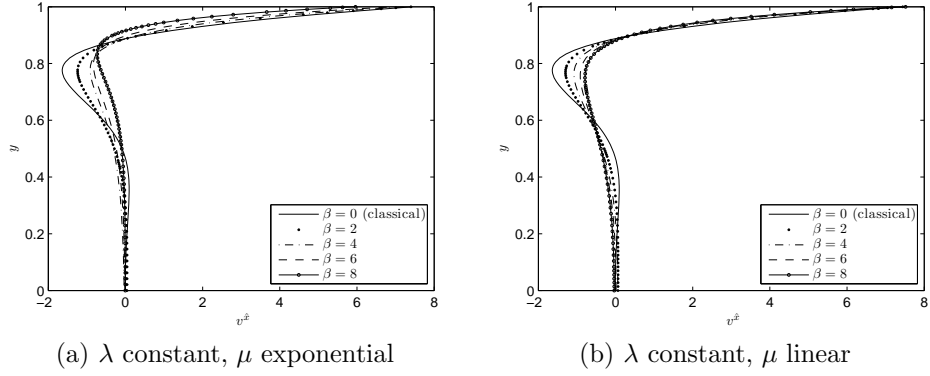

(a) $\lambda$ constant, $\mu$ exponential

(b) $\lambda$ constant, $\mu$ linear

Figure 9: Vorticity profiles at $\omega t = \frac{\pi}{2}$ for various $\beta$ values when the viscosity ($\mu$) depends on pressure exponentially and linearly. Relaxation time ($\lambda$) is kept constant in both cases with Re = 100, We = 0.1, $\omega = 1$, $\Pi = 1$.

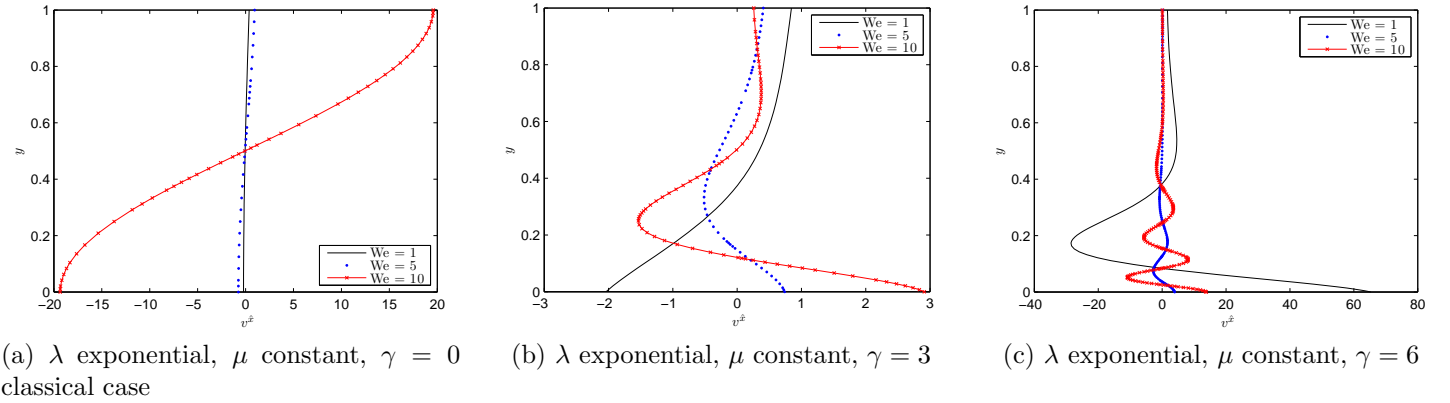

(a) $\lambda$ exponential, $\mu$ constant, $\gamma = 0$ classical case

(b) $\lambda$ exponential, $\mu$ constant, $\gamma = 3$

(c) $\lambda$ exponential, $\mu$ constant, $\gamma = 6$

Figure 10: Vorticity profiles at $\omega t = \frac{\pi}{2}$ for various We values for the classical case as well as when the relaxation time ($\lambda$) depends on pressure exponentially. Viscosity ($\mu$) is kept constant in all cases with Re = 1, $\omega = 1$, $\Pi = 1$.

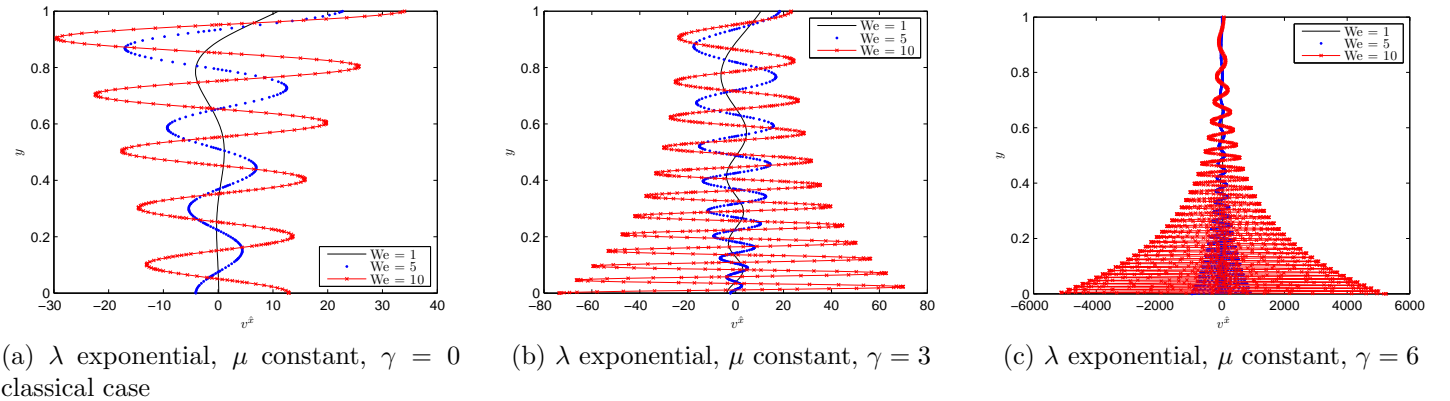

(a) $\lambda$ exponential, $\mu$ constant, $\gamma = 0$ classical case

(b) $\lambda$ exponential, $\mu$ constant, $\gamma = 3$

(c) $\lambda$ exponential, $\mu$ constant, $\gamma = 6$

Figure 11: Vorticity profiles at $\omega t = \frac{\pi}{2}$ for various We values for the classical case as well as when the relaxation time ($\lambda$) depends on pressure exponentially. Viscosity ($\mu$) is kept constant in all cases with Re = 100, $\omega = 1$, $\Pi = 1$.

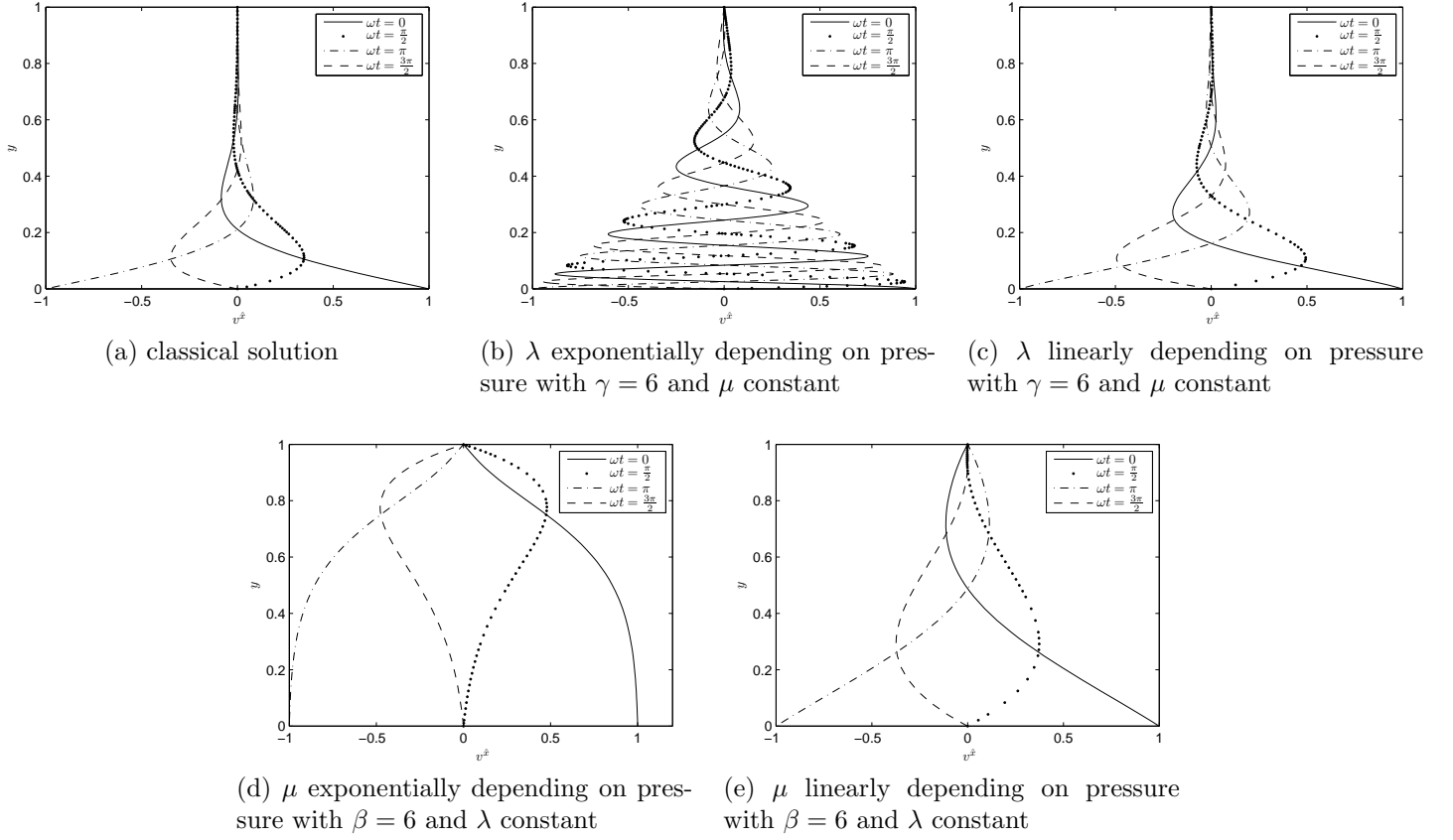

Figure 12: Snapshots of velocity profiles for various cases of relaxation time ($\lambda$) and viscosity ($\mu$) depending on pressure, with bottom plate oscillating and for Re = 100, We = 0.1, $\omega = 1$, $\Pi = 1$.

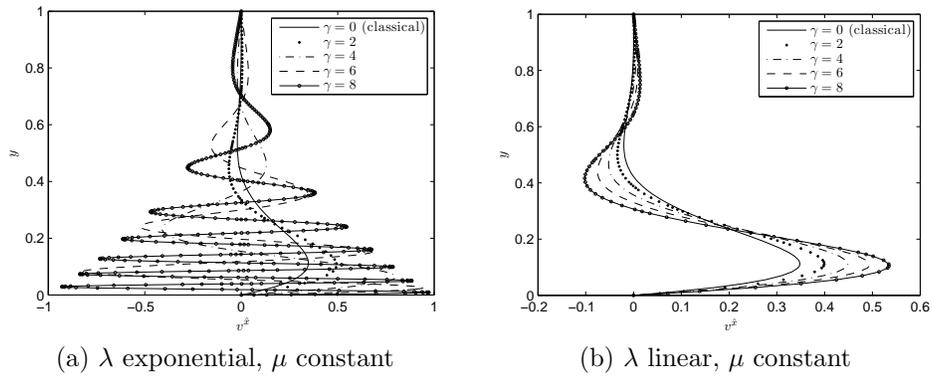

Figure 13: Velocity profiles at $\omega t = \frac{\pi}{2}$ for various $\gamma$ values when the relaxation time ($\lambda$) depends on pressure exponentially and linearly, when the bottom plate is oscillating. Viscosity ($\mu$) is kept constant in both cases with Re = 100, We = 0.1, $\omega = 1$, $\Pi = 1$.

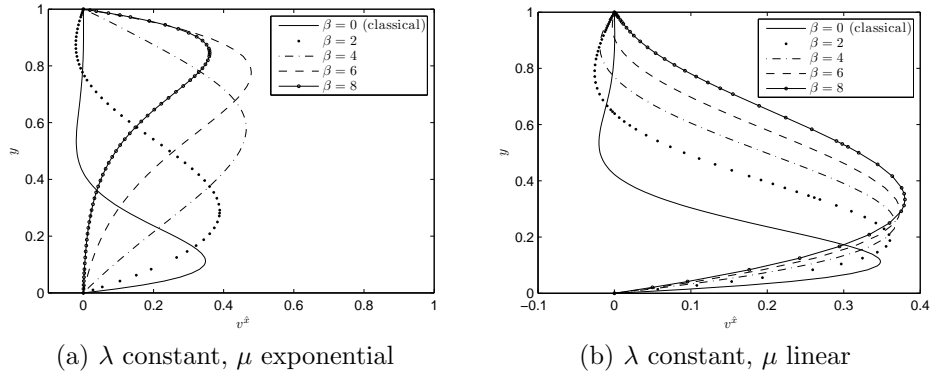

(a) $\lambda$ constant, $\mu$ exponential
(b) $\lambda$ constant, $\mu$ linear

Figure 14: Velocity profiles at $\omega t = \frac{\pi}{2}$ for various $\beta$ values when the viscosity ($\mu$) depends on pressure exponentially and linearly, when the bottom plate is oscillating. Relaxation time ($\lambda$) is kept constant in both cases with Re = 100, We = 0.1, $\omega = 1$, $\Pi = 1$.

can be expected if the fact that the material moduli depend on pressure is taken into account in studying the motion of viscoelastic fluids.

# References

[1] J. C. Maxwell, On the Dynamical Theory of Gases, Philos. Trans. R. Soc. 157 (1867) 49–88, doi: 10.1098/rstl.1867.0004.

[2] J. D. Ferry, Viscoelastic properties of polymers, John Wiley & Sons, New York, 3 edn., 1980.

[3] L. M. Cathles, The Viscosity of the Earth's mantle, Princeton University Press, Princeton, 1975.

[4] P. W. Bridgman, The physics of high pressure, Macmillan, New York, 1931.

[5] M. J. Neale (Ed.), Tribology handbook, John Wiley & Sons, Inc., New York, 1973.

[6] D. Gwynllyw, A. Davies, T. Phillips, On the effects of a piezoviscous lubricant on the dynamics of a journal bearing, J. Rheol. 40 (6) (1996) 1239–1266, ISSN 0148-6055.

[7] H. Singh, A. W. Nolle, Pressure dependence of the viscoelastic behavior of polyisobutylene, J. Appl. Phys. 30 (3) (1959) 337–341, ISSN 0021-8979.

[8] J. E. McKinney, H. V. Belcher, Dynamic compressibility of polyvinylacetate and its relation to free volume, J. Res. Nat. Bur. Stand. Sect. A. Phys. Chem. A 67 (1) (1963) 43–53.

[9] J. Weertman, S. White, A. H. Cook, Creep Laws for the Mantle of the Earth [and Discussion], Proc. R. Soc. Lond., Ser. A, Math. Phys. Eng. Sci. 288 (1350) (1978) 9–26, ISSN 00804614, URL http://www.jstor.org/stable/74973.

[10] E. R. Ivins, C. G. Sammis, C. F. Yoder, Deep mantle viscous structure with prior estimate and satellite constraint, J. Geophys. Res. 98 (1993) 4579–4609.

[11] T. Sahaphol, S. Miura, Shear moduli of volcanic soils, Soil Dyn. Earthq. Eng. 25 (2) (2005) 157–165, ISSN 0267-7261, doi:10.1016/j.soildyn.2004.10.001, URL http://www.sciencedirect.com/science/article/B6V4Y-4F0GR6P-1/2/57e03f9fbd55ce70669a15555fe4

[12] R. K. McConnel, Isostatic adjustment in a layered earth, J. Geophys. Res. 70 (20) (1965) 5171–5188, ISSN 0148-0227.

[13] R. K. McConnel, Viscosity of mantle from relaxation time spectra of isostatic adjustment, J. Geophys. Res. 73 (22) (1968) 7089–7105, ISSN 0148-0227.

[14] P. Wu, H. Wang, Postglacial isostatic adjustment in a self-gravitating spherical Earth with power-law rheology, J. Geodyn. 46 (3-5) (2008) 118–130, ISSN 0264-3707, doi:10.1016/j.jog.2008.03.008.

[15] K. R. Rajagopal, A. R. Srinivasa, On the thermodynamics of fluids defined by implicit constitutive relations, Z. Angew. Math. Phys. 59 (4) (2008) 715–729, doi:10.1007/s00033-007-7039-1.

[16] K. R. Rajagopal, On implicit constitutive theories for fluids, J. Fluid Mech. 550 (2006) 243–249, doi:10.1017/S0022112005008025.

[17] J. Hron, J. Málek, K. R. Rajagopal, Simple flows of fluids with pressure-dependent viscosities, Proc. R. Soc. Lond., Ser. A, Math. Phys. Eng. Sci. 457 (2011) (2001) 1603–1622.

[18] S. Srinivasan, K. R. Rajagopal, Study of a variant of Stokes' first and second problems for fluids with pressure dependent viscosities, Int. J. Eng. Sci. 47 (11-12) (2009) 1357–1366, ISSN 0020-7225, doi:10.1016/j.ijengsci.2008.11.002, URL http://www.sciencedirect.com/science/article/B6V32-4V6RNRG-1/2/562d1ab4be6cd8a91083a5ec28cf

[19] K. R. Rajagopal, A. R. Srinivasa, A thermodynamic frame work for rate type fluid models, J. Non-Newton. Fluid Mech. 88 (3) (2000) 207–227, ISSN 0377-0257, doi:10.1016/S0377-0257(99)00023-3, URL http://www.sciencedirect.com/science/article/B6TGV-3Y9RCSK-1/2/c9f229a12a7f3f888fa3f460c58e

[20] G. G. Stokes, On the effect of the internal friction of fluids on the motion of pendulums, Trans. Cambridge Phil. Soc. 9 (1851) 8–106.

[21] N. A. Haskell, The motion of a viscous fluid under a surface load, Physics-J. Gen. Appl. Phys. 6 (1) (1935) 265–269, ISSN 0148-6349.

[22] A. E. H. Love, Some problems of geodynamics, Cambridge University Press, Cambridge, 1911.

[23] D. Wolf, Lamé's problem of gravitational viscoelasticity: the isochemical, incompressible planet, Geophys. J. Int. 117 (2) (1994) 321–348, doi:10.1111/j.1365-246X.1994.tb01801.x.

[24] C. Fetecau, C. Fetecau, A new exact solution for the flow of a Maxwell fluid past an infinite plate, Int. J. Non-Linear Mech. 38 (3) (2003) 423–427, ISSN 0020-7462, doi:10.1016/S0020-7462(01)00062-2.

[25] A. D. Polyanin, V. F. Zaitsev, Handbook of exact solutions for ordinary differential equations, Chapman & Hall/CRC, Boca Raton, FL, second edn., ISBN 1-58488-297-2, 2003.

[26] G. N. Watson, A treatise on the theory of Bessel functions, Cambridge University Press, 2 edn., reprinted 1980, 1944.

[27] D. E. Amos, Algorithm 644: A Portable Package for Bessel Functions of a Complex Argument and Nonnegative Order, ACM Trans. Math. Softw. 12 (3) (1986) 265–273, ISSN 0098-3500, URL http://doi.acm.org/10.1145/7921.214331.

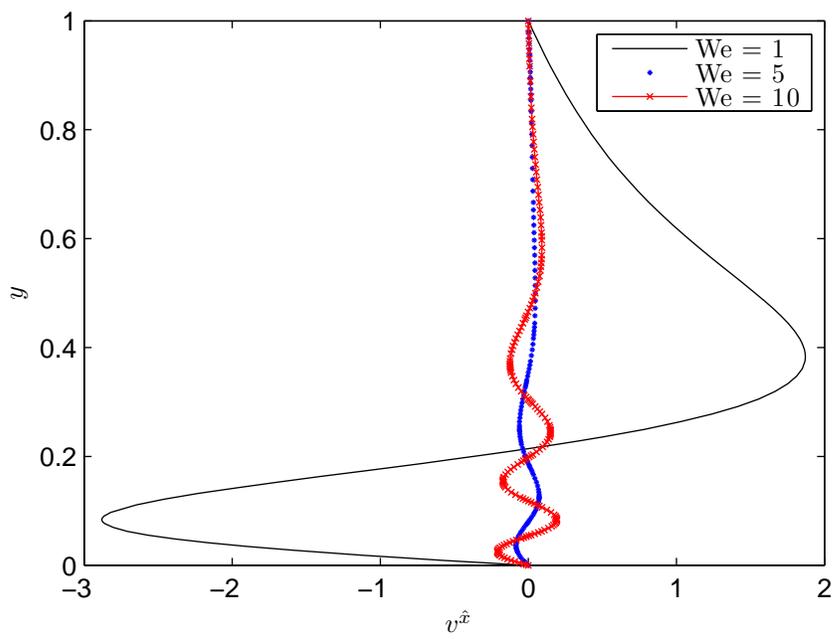